\documentclass[12pt]{article}
\usepackage{color}
\definecolor{darkblue}{rgb}{0.00,0.25,0.50}
\usepackage[colorlinks,filecolor=blue,citecolor=darkblue]{hyperref}

\usepackage{amsfonts,amssymb,amsmath,amsthm}
\usepackage{url}
\usepackage{enumerate}
\usepackage[ukrainian, russian, english]{babel}
\usepackage[cp1251]{inputenc}
\begin{document}\selectlanguage{ukrainian}
\thispagestyle{empty}

\title{}

\begin{center}
\textbf{\Large Оцінки найкращих наближень  класів нескінченно
диференційовних функцій в рівномірній та інтегральних метриках}
\end{center}
\vskip0.5cm
\begin{center}
А.~С.~Сердюк${}^1$, Т.~А.~Степанюк${}^2$\\ \emph{\small
${}^1$Інститут математики НАН
України, Київ\\
${}^2$Східноєвропейський національний університет імені Лесі
Українки, Луцьк\\}
\end{center}
\vskip0.5cm

%\address{Institute of Mathematics of the National Academy of
%Sciences of Ukraine\\ 3\\ Tereshenkivska st.\\ 01601\\ Kiev, Ukraine}

\begin{abstract}
Найдены равномерные относительно  параметра $p \ (1\leq p\leq\infty)$  оценки сверху наилучших приближений
тригонометрическими полиномами классов периодических функций
$C^{\psi}_{\beta,p}$, порождаемых последовательностями $\psi(k)$,  убывающих к нулю быстрее  любой
степенной функции.
Полученные оценки точны  по порядку и содержат выраженные в явном виде постоянные,  зависящие только от функции
$\psi$.
 Аналогичные оценки установлены для наилучших приближений  классов
$L^{\psi}_{\beta,1}$
в метриках пространств $L_{s}$, $1\leq s\leq \infty$.

\vskip 0.5cm

 We find uniform with respect to parameter  $p \ (1\leq p\leq\infty)$ upper estimations of best approximations by trigonometric polynomials  of classes $C^{\psi}_{\beta,p}$  of periodic functions generated by
sequences   $\psi(k)$, that decrease to nought
faster than any power function. Obtained estimations are exact for order and have constants that are written in explicit form and
depend on function $\psi$ only.
We obtain analogical estimations for best approximations of classes
$L^{\psi}_{\beta,1}$
in metrics of spaces  $L_{s}$, ${1\leq s\leq \infty}$.

\end{abstract}

%%%%%%%%%%%%%%%%%%%%%%%%%%%%%%%%%%%%%%%%%%%%%%%%%%%%%%%%%%%%%%%%%%%%%%%%%

Нехай $C$ --- простір $2\pi$--періодичних неперервних функцій, у
якому норма задається за допомогою рівності $
{\|f\|_{C}=\max\limits_{t}|f(t)|}$; $L_{\infty}$ --- простір
$2\pi$--періодичних вимірних і суттєво обмежених функцій $f(t)$ з нормою
$\|f\|_{\infty}=\mathop{\rm{ess}\sup}\limits_{t}|f(t)|$; $L_{p}$,
$1\leq p<\infty$, --- простір $2\pi$--періодичних сумовних в $p$--му
степені на $[0,2\pi)$ функцій $f(t)$, в якому  норма задана формулою
${\|f\|_{p}=\Big(\int\limits_{0}^{2\pi}|f(t)|^{p}dt\Big)^{\frac{1}{p}}}.$

Нехай, далі, $L^{\psi}_{\beta,p}, \ 1\leq p\leq \infty$,
---
 множина всіх $2\pi$--періодичних функцій $f$, котрі майже для всіх
$x\in\mathbb{R}$ зображуються за допомогою згортки
\begin{equation}\label{zgo}
f(x)=\frac{a_{0}}{2}+\frac{1}{\pi}\int\limits_{-\pi}^{\pi}\Psi_{\beta}(x-t)\varphi(t)dt,
\ a_{0}\in\mathbb{R},  \ \varphi\in B^{0}_{p},
\end{equation}
де  $\Psi_{\beta}(t)$
--- сумовна на $[0,2\pi)$ функція, ряд Фур'є якої має вигляд
$$
\sum\limits_{k=1}^{\infty}\psi(k)\cos
\big(kt-\frac{\beta\pi}{2}\big), \ \  \beta\in
    \mathbb{R}, \ \psi(k)>0,
$$
а
$$B_{p}^{0}=\left\{\varphi\in L_{1}: \ ||\varphi||_{p}\leq 1, \
\varphi\perp 1\right\}, \ \ 1\leq p\leq \infty.$$

 Функцію
$\varphi$ в зображенні (\ref{zgo}), згідно з О.І. Степанцем \cite[с. 132]{Stepanets1}, називають $(\psi,\beta)$--похідною
функції $f$ і позначають через $f^{\psi}_{\beta}$.
Підмножину неперервних функцій із $L^{\psi}_{\beta,p}, \ {1\leq p\leq
\infty,}$ позначають через $C^{\psi}_{\beta,p}, \ 1\leq p\leq
\infty$.

Послідовності $\psi (k),\ k\in \mathbb{N},$ що визначають класи
$L^{\psi}_{\beta,p}$ та $C^{\psi}_{\beta,p}$, зручно розглядати як звуження на множину
натуральних чисел $\mathbb{N}$ деяких додатних, неперервних, опуклих
донизу функцій $\psi(t)$ неперервного аргументу $t\geq1$ таких, що $
\lim\limits_{t\rightarrow\infty}\psi(t)=0. $ Множину всіх таких
функцій $\psi(t)$ позначатимемо через ${\mathfrak M}$.

Наслідуючи О.І. Степанця (див., наприклад, \cite[с.
160]{Stepanets1}), кожній функції $\psi\in{\mathfrak M}$ поставимо у відповідність
характеристики
$$
\eta(t)=\eta(\psi;t)=\psi^{-1}\left(\psi(t)/2\right), \ \ \
\mu(t)=\mu(\psi;t)=\frac{t}{\eta(t)-t},
$$
де $\psi^{-1}$ --- обернена до $\psi$ функція i покладемо
$$
\mathfrak{M}^{+}_{\infty}=\left\{\psi\in \mathfrak{M}: \ \
\mu(\psi;t)\uparrow\infty, \ t\rightarrow\infty \right\}.
$$

Якщо $\psi\in\mathfrak{M}^{+}_{\infty}$, то (див., наприклад,
\cite[с. 97]{Step monog 1987}) функція $\psi(t)$ спадає до нуля швидше довільної степеневої функції, тобто:
$$
\forall r\in\mathbb{R} \ \ \lim\limits_{t\rightarrow\infty}t^{r}\psi(t)=0.
$$
Це означає, що за умови $\psi\in\mathfrak{M}^{+}_{\infty}$, ряд Фур'є довільної функції $f$ із
 $C^{\psi}_{\beta,p}, \ \beta\in\mathbb{R}$ можна диференціювати довільне число разів і в результаті будемо одержувати рівномірно збіжні ряди. Отже, класи  $C^{\psi}_{\beta,p}$ при $\psi\in\mathfrak{M}^{+}_{\infty}$
складаються з нескінченно диференційовних функцій.

З іншого боку, як
показано в \cite[с. 1692]{Stepanets_Serdyuk_Shydlich}, для кожної
нескінченно диференційовної, $2\pi$--періодичної функції $f$ можна
вказати функцію $\psi$ з множини $\mathfrak{M}^{+}_{\infty}$ таку,
що $f\in C^{\psi}_{\beta,p}$ для довільних $\beta\in\mathbb{R}$.

Із множини  $\mathfrak{M}^{+}_{\infty}$  виділяють підмножини ${\mathfrak M^{'}_{\infty}}$ і ${\mathfrak M^{''}_{\infty}}$; ${\mathfrak M^{'}_{\infty}}$ --- множина функцій
${\psi\in {\mathfrak M^{+}_{\infty}}}$, для яких величина
$\eta(\psi;t)-t$ обмежена зверху, тобто існує стала $K_{1}>0$ така,
що ${\eta(\psi;t)-t\leq K_{1}}, \ t\geq1$, а ${\mathfrak
M^{''}_{\infty}}$ --- множина функцій ${\psi\in {\mathfrak
M^{+}_{\infty}}}$, для яких величина $\eta(\psi;t)-t$ обмежена знизу
деяким додатним числом, тобто існує стала $K_{2}>0$ така, що
${\eta(\psi;t)-t\geq K_{2}}, \ t\geq1$.

Типовими представниками множини $\mathfrak{M}^{+}_{\infty}$  є
функції ${\psi_{r, \alpha}(t)=\exp(-\alpha t^{r})}, \ {\alpha>0}, \ r>0$,
причому, якщо $r\geq1$, то $\psi_{r,\alpha}\in \mathfrak{M}^{'}_{\infty}$,
а якщо ${r\in(0,1]}$, то $\psi_{r,\alpha}\in \mathfrak{M}^{''}_{\infty}$.
Класи $C^{\psi}_{\beta,p}$ та $L^{\psi}_{\beta,p}$, що породжуються функціями $\psi=\psi_{r, \alpha}$
будемо позначати через $C^{\alpha,r}_{\beta,p}$ та $L^{\alpha,r}_{\beta,p}$ відповідно.

Нехай, далі, ${E}_{n}(C^{\psi}_{\beta,p})_{C}$ та ${E}_{n}(L^{\psi}_{\beta,p})_{s}$ --- найкращі наближення класів
$C^{\psi}_{\beta,p}$ та $L^{\psi}_{\beta,p}$ в метриках просторів $C$ та $L_{s}$, тобто величини вигляду
$$
{E}_{n}(C^{\psi}_{\beta,p})_{C}=\sup\limits_{f\in
C^{\psi}_{\beta,p}}\inf\limits_{t_{n-1}\in\mathcal{T}_{2n-1}}\|f(\cdot)-t_{n-1}(\cdot)\|_{C},
\ 1\leq p\leq \infty,
$$
$$
{E}_{n}(L^{\psi}_{\beta,p})_{s}=\sup\limits_{f\in
L^{\psi}_{\beta,p}}\inf\limits_{t_{n-1}\in\mathcal{T}_{2n-1}}\|f(\cdot)-t_{n-1}(\cdot)\|_{s},
\ 1\leq p,s\leq \infty,
$$
де $\mathcal{T}_{2n-1}$ --- підпростір усіх тригонометричних
поліномів $t_{n-1}$ порядку не вищого за ${n-1}$.
 Дана робота присвячена
знаходженню точних порядкових оцінок величин ${E}_{n}(C^{\psi}_{\beta,p})_{C}$ i ${E}_{n}(L^{\psi}_{\beta,1})_{s}$ при
$\psi\in\mathfrak{M^{+}_{\infty}}$ i $\beta\in\mathbb{R}$.

В роботі  \cite[с. 225]{Step monog 1987} (див. також \cite[с. 48]{Stepanets2}) встановлено  точні порядкові
оцінки величин $E_{n}(L^{\psi}_{\beta,p})_{s}$, $\beta\in\mathbb{R}$, при $\psi\in {\mathfrak M^{'}_{\infty}}$, $1\leq
p,s\leq\infty$, а також
при $\psi\in {\mathfrak M^{''}_{\infty}}$,
${1<p,s<\infty}$.
У випадку $p=s=1$ або  ${p=s=\infty}$ в \cite{serdyuk2004zbirnyk}
встановлено асимптотичні рівності при $n\rightarrow\infty$ для величин
$E_{n}(L^{\psi}_{\beta,p})_{s}$,  ${\psi\in {\mathfrak
M^{+}_{\infty}}}$ і $\beta\in\mathbb{R}$.
Якщо ж для послідовності $\psi(k), \ {k\in\mathbb{N}}$
виконуються умови: 1)~${\Delta^{2}\psi(k)\mathop{=}\limits^{\rm
df}\psi(k)-2\psi(k+1)+\psi(k+2)\geq 0}$,
${\frac{\psi(k+1)}{\psi(k)}\leq \rho}, \ {0<\rho<1}, \ k=n,
n+1,...$; 2)~${\frac{\Delta^{2}\psi(n)}{\psi(n)}>\frac{(1+3\rho)\rho^{2n}}{(1-\rho)\sqrt{1-2\rho^{2n}}}
} $, то при $p=s=1$,  $p=s=\infty$, в \cite{Serdyuk2002}  отримано
точні значення величин $E_{n}(L^{\psi}_{\beta,p})_{s}$  для довільних $\beta\in\mathbb{R}$
і $n\in\mathbb{N}$.

В роботі
\cite{S_S} знайдено порядкові оцінки величин
${E}_{n}(C^{\psi}_{\beta,p})_{C}, \ 1\leq p<\infty,$ і
${E}_{n}(L^{\psi}_{\beta,1})_{s}, \  {1<s\leq \infty},$ у випадку,
коли $\psi\in\mathfrak{M}^{+}_{\infty}$,
${\eta(t)-t\geq a>2}, \ {\mu(t)\geq b>2}$.
При цьому для оцінок зверху для найкращих наближень  ${E}_{n}(C^{\psi}_{\beta,p})_{C}$ i ${E}_{n}(L^{\psi}_{\beta,1})_{s}$ використовувались відповідно величини
 $$
  {\cal E}_{n}(C^{\psi}_{\beta,p})_{C}=\sup\limits_{f\in
C^{\psi}_{\beta,p}}\|f(\cdot)-S_{n-1}(f;\cdot)\|_{C},  \ \ 1\leq p<\infty,
  $$
  $$
  {\cal E}_{n}(L^{\psi}_{\beta,1})_{s}=\sup\limits_{f\in
L^{\psi}_{\beta,1}}\|f(\cdot)-S_{n-1}(f;\cdot)\|_{s}, \ \ 1<s\leq\infty,
  $$
де $S_{n-1}(f;\cdot)$ --- частинні суми Фур'є порядку $n-1$ функції $f$.
Було показано справедливість нерівностей
\begin{equation}\label{th1}
{\cal
E}_{n}(C^{\psi}_{\beta,p})_{C} \leq C_{a,b} \
(2p)^{1-\frac{1}{p}} \psi(n)(\eta(n)-n)^{\frac{1}{p}}, \ \ 1\leq p<\infty,
\end{equation}
\begin{equation}\label{th2}
{\cal
E}_{n}(L^{\psi}_{\beta,1})_{s}  \leq C_{a,b} \
\left(2s'\right)^{\frac{1}{s}}\psi(n)\left(\eta(n)-n\right)^{\frac{1}{s'}}, \ \ 1<s\leq\infty,
\  \frac{1}{s}+\frac{1}{s'}=1,
\end{equation}
де
\begin{equation}\label{Cab_F}
C_{a,b}=
\frac{1}{\pi}\max\Big\{\frac{2b}{b-2}+\frac{1}{a}, \ 2\pi\Big\}.
\end{equation}
Крім того, в \cite{S_S} доведено, що при
 $\psi\in\mathfrak{M}^{+}_{\infty}$,
${\eta(n)-n>2}$  суми Фур'є забезпечують порядок найкращих наближень тригонометричними поліномами, тобто
$$
  { E}_{n}(C^{\psi}_{\beta,p})_{C}\asymp{\cal E}_{n}(C^{\psi}_{\beta,p})_{C}\asymp\psi(n)(\eta(n)-n)^{\frac{1}{p}},
  \ \ 1\leq p<\infty,
  $$
  $$
  { E}_{n}(L^{\psi}_{\beta,1})_{s}\asymp{\cal E}_{n}(L^{\psi}_{\beta,1})_{s}\asymp\psi(n)(\eta(n)-n)^{\frac{1}{s'}}, \ \ 1<s\leq\infty, \ \frac{1}{s}+\frac{1}{s'}=1,
  $$
  (тут і надалі запис $A(n)\asymp B(n)$ $(A(n)>0, \ B(n)>0)$ означає існування
додатних сталих $K_{1}$ і $K_{2}$ таких, що ${K_{1}B(n)\leq A(n)\leq
K_{2}B(n)}$, $n\in\mathbb{N}$).

 Зазначимо, що порядкова рівність ${\cal E}_{n}(C^{\psi}_{\beta,p})_{C}\asymp\psi(n)(\eta(n)-n)^{\frac{1}{p}},
 \ 1\leq p<\infty,$ у випадку $\psi\in\mathfrak{M}_{\infty}^{''}$ встановлена в роботі  \cite{Rom}.
 Крім того, у роботах \cite{Serdyuk} і \cite{Serdyuk2013}
 при $p=2$, $s=\infty$ та $p=1$, $s=2$ за умови збіжності ряду
$\sum\limits_{k=1}^{\infty}\psi^{2}(k)$ знайдено точні значення величин
${\cal E}_{n}(L^{\psi}_{\beta,p})_{s}$  для всіх $\beta\in\mathbb{R}$ і $n\in\mathbb{N}$.

Водночас, суми Фур'є, як апарат наближення, не дозволяють записати рівномірних відносно параметрів $p, \ 1\leq p\leq\infty, $ та $s, \ 1\leq s\leq\infty,$ оцінок зверху для величин ${
E}_{n}(C^{\psi}_{\beta,p})_{C}$ та ${
E}_{n}(L^{\psi}_{\beta,1})_{s}$. Цей факт зумовлений тією обставиною, що при $p=\infty$ та $s=1$ за умови
$\psi\in\mathfrak{M}_{\infty}^{''}$
$$
\frac{{\cal E}_{n}(C^{\psi}_{\beta,p})_{C}}{{ E}_{n}(C^{\psi}_{\beta,p})_{C}}\asymp\ln^{+}(\eta(n)-n), \ \ \ \ \
\frac{{\cal E}_{n}(L^{\psi}_{\beta,1})_{s}}{{ E}_{n}(L^{\psi}_{\beta,1})_{s}}\asymp\ln^{+}(\eta(n)-n),
$$
де $\ln^{+}t=\max\{0, \ \ln t\}$ (див., наприклад, \cite[с. 264, 285]{Stepanets1} та
\cite[с. 86,87]{Stepanets2}.

В даній роботі побудовано лінійний метод наближення $V_{n,\psi}(t)$, що дозволяє записати рівномірні відносно параметрів $p \ (1\leq p\leq\infty)$ i $s \ (1\leq s\leq\infty)$ оцінки зверху найкращих наближень
${E}_{n}(C^{\psi}_{\beta,p})_{C}$ i ${E}_{n}(L^{\psi}_{\beta,1})_{s}$ при $\psi\in {\mathfrak
M^{+}_{\infty}}$ і $\beta\in\mathbb{R}$.
Показано, що знайдені оцінки є точними за порядком і
містять виражені в явному вигляді сталі, які залежать лише від функції
$\psi$.

\textbf{Теорема 1.} \emph{ Нехай $\psi\in\mathfrak{M}^{+}_{\infty}$,
$\beta\in \mathbb{R}$, ${1\leq p\leq \infty}$. Тоді  для
$n\in
 \mathbb{N}$ таких, що $ {\eta(n)-n\geq a>2}, \ {\mu(n)\geq b>2}$ справедливі
 оцінки
\begin{equation}\label{theorem1}
C_{a}\psi(n)(\eta(n)-n)^{\frac{1}{p}}\leq{
E}_{n}(C^{\psi}_{\beta,p})_{C}\leq
C^{*}_{a,b}\psi(n)(\eta(n)-n)^{\frac{1}{p}},
\end{equation}
де}
\begin{equation}\label{Ca}
C_{a}=\frac{\pi}{96\left(1+\pi^{2}\right)^{2}}
\frac{(a-1)^{2}(a-2)^{2}}{a^{3}(3a-4)},
\end{equation}
\begin{equation}\label{Cab}
C^{*}_{a,b}=\frac{2(1+\pi^{2})}{\pi}\Big(\frac{2b}{b-2}+\frac{a}{a-1}\Big).
\end{equation}

\textbf{\emph{Доведення теореми 1.}} Згідно з формулою (33)
роботи \cite{S_S} для довільної ${\psi\in\mathfrak{M}^{+}_{\infty}}$ при $ {\eta(n)-n\geq a>2}, \ {\mu(n)\geq b>2}$
має місце оцінка
$$
{
E}_{n}(C^{\psi}_{\beta,p})_{C}\geq C_{a}\psi(n)(\eta(n)-n)^{\frac{1}{p}}, \ \ 1\leq p<\infty,
$$
в якій величина $C_{a}$ означена рівністю(\ref{Ca}). Для знаходження оцінки зверху величини
${
E}_{n}(C^{\psi}_{\beta,p})_{C}$
розглянемо
лінійний метод $V_{n,\psi}(f;x)$ наближення функцій з множини $C^{\psi}_{\beta,p}(L^{\psi}_{\beta,1})$ наступного  вигляду:
\begin{equation}\label{sum}
V_{n,\psi}(f;x)=\frac{a_{0}(f)}{2}+\sum\limits_{k=1}^{n-1}\lambda_{n,[\eta(n)]-n+1}(k)\left(a_{k}(f)\cos
kx+b_{k}(f)\sin kx\right), \ n\in\mathbb{N},
\end{equation}
де
\begin{equation}\label{lamb}
\lambda_{n,[\eta(n)]-n+1}(k)={\left\{\begin{array}{cc}
1, \ & 0 \leq k \leq 2n-[\eta(n)]-1, \\
1-\frac{[\eta(n)]-2n+k}{[\eta(n)]-n}\frac{\psi(n)}{\psi(k)}, &
2n-[\eta(n)] \leq  k\leq n-1 \
  \end{array} \right.}
\end{equation}
(тут і надалі $[\alpha]$ --- ціла частина дійсного числа $\alpha$).

Зауважимо, що суми (\ref{sum}) є частинним випадком узагальнених сум
Валле Пуссена $U_{n,m}^{\nu}(f;x)$, тобто поліномів вигляду (див., наприклад,  \cite{Serdyuk_ovsii}:
$$
U_{n,m}^{\nu}(f;x)=\frac{a_{0}(f)}{2}+\sum\limits_{k=1}^{n-1}\lambda_{n,m}(k)\left(a_{k}(f)\cos
kx+b_{k}(f)\sin kx\right),
$$
де
$$
 \lambda_{n,m}(k)={\left\{\begin{array}{cc}
1, \ &  0\leq k\leq n-m, \\
1-\frac{\nu(k)}{\nu(n)}, & n-m+1\leq k\leq n-1, \
  \end{array} \right.}
$$
а $\nu(k), \ k\in\mathbb{N}$ --- довільно монотонно зростаюча послідовність додатних
 чисел ${m\in[1,n]}, \ {m\in\mathbb{N}}$. Поклавши
$\nu(k)=\frac{[\eta(n)]-2n+k}{\psi(k)}$ і $m=[\eta(n)]-n+1$, одержимо
${U_{n,m}^{\nu}(f;x)=V_{n,\psi}(f;x)}$.

Розглянемо величину
$$
{\cal
E}(C^{\psi}_{\beta,p};V_{n,\psi})_{C}=\sup\limits_{f\in
C^{\psi}_{\beta,p}}\|f(\cdot)-V_{n,\psi}(f;\cdot)\|_{C}, \ \ 1\leq
p\leq \infty.
$$
Оскільки
\begin{equation}\label{dop}
{ E}_{n}(C^{\psi}_{\beta,p})_{C}\leq {\cal
E}(C^{\psi}_{\beta,p};V_{n,\psi})_{C},
\end{equation}
то для доведення теореми 1 достатньо показати, що при $\psi\in\mathfrak{M}^{+}_{\infty}$, $ {\eta(n)-n\geq a>2}, \ {\mu(n)\geq b>2}$ має місце
співвідношення
\begin{equation}\label{dop1}
{\cal E}(C^{\psi}_{\beta,p};V_{n,\psi})_{C}\leq
C_{a,b}^{*}\psi(n)(\eta(n)-n)^{\frac{1}{p}}, \ 1\leq p\leq \infty,
\end{equation}
де $C_{a,b}^{*}$ означаються рівністю (\ref{Cab}).
Згідно з \cite[с. 51]{Step monog 1987} для усіх $x\in\mathbb{R}$ має місце
рівність
$$
f(x)-V_{n,\psi}(f;x)=\frac{1}{\pi}\int\limits_{-\pi}^{\pi}f^{\psi}_{\beta}(x+t)
\left(\sum\limits_{k=1}^{n-1}(1-\lambda_{n,[\eta(n)]-n+1}(k))\psi(k)\cos\Big(kt+\frac{\beta\pi}{2}\Big)+\right.
$$
\begin{equation}\label{rizn}
\left.+\sum\limits_{k=n}^{\infty}\psi(k)\cos\Big(kt+\frac{\beta\pi}{2}\Big)\right)dt,
\end{equation}
де $\lambda_{n,[\eta(n)]-n+1}(k)$  означаються  формулою
(\ref{lamb}).  Із (\ref{rizn}) і (\ref{lamb}), отримаємо
\begin{equation}\label{fff}
f(x)-V_{n,\psi}(f;x)=\frac{1}{\pi}\int\limits_{-\pi}^{\pi}f^{\psi}_{\beta}(t)\Psi^{*}_{\beta,n}(x-t)dt, \ \|f^{\psi}_{\beta}\|_{p}\leq1, \  1\leq p\leq\infty,
\end{equation}
де
$$
\Psi^{*}_{\beta,n}(t)=
\psi(n)\sum\limits_{k=2n-[\eta(n)]}^{n-1}\frac{[\eta(n)]-2n+k}{[\eta(n)]-n}\cos\Big(kt-\frac{\beta\pi}{2}\Big)
+\sum\limits_{k=n}^{\infty}\psi(k)\cos\Big(kt-\frac{\beta\pi}{2}\Big)=
$$
\begin{equation}\label{kernel}
=\psi(n)\sum\limits_{k=2n-[\eta(n)]+1}^{n-1}\Big(1-\frac{n-k}{[\eta(n)]-n}\Big)\cos\Big(kt-\frac{\beta\pi}{2}\Big)
+\sum\limits_{k=n}^{\infty}\psi(k)\cos\Big(kt-\frac{\beta\pi}{2}\Big).
\end{equation}
Згідно з  твердженнями 8.1 і 8.2 роботи \cite[с. 137--138]{Stepanets1}  маємо, що для будь-якої функції $f\in C_{\beta,p}^{\psi}$ виконується нерівність
\begin{equation}\label{zz2}
\Bigg\|\frac{1}{\pi}\int\limits_{-\pi}^{\pi}f^{\psi}_{\beta}(t)\Psi^{*}_{\beta,n}(x-t)dt\Bigg\|_{C}
\leq\frac{1}{\pi}\|f^{\psi}_{\beta}\|_{p}\|\Psi^{*}_{\beta,n}\|_{p'}\leq
\frac{1}{\pi}\|\Psi^{*}_{\beta,n}(t)\|_{p'}, \ \ 1\leq
p\leq \infty, \ \ \frac{1}{p}+\frac{1}{p'}=1.
\end{equation}

Оцінимо зверху величину $\|\Psi^{*}_{\beta,n}\|_{p'}$.
Позначивши
\begin{equation}\label{Dkb}
D_{k,\beta}(t)=\frac{1}{2}\cos\frac{\beta\pi}{2}+\sum\limits_{j=1}^{k}\cos\Big(jt-\frac{\beta\pi}{2}\Big),
\end{equation}
з (\ref{kernel}) одержуємо
$$
\Psi^{*}_{\beta,n}(t)=
\psi(n)\left(D_{n-1,\beta}(t)-D_{2n-[\eta(n)],\beta}(t)\right)-
\frac{\psi(n)}{[\eta(n)]-n}\sum\limits_{k=2n-[\eta(n)]+1}^{n-1}(n-k)\cos\Big(kt-\frac{\beta\pi}{2}\Big)+
$$
\begin{equation}\label{q}
+\sum\limits_{k=n}^{\infty}\psi(k)\cos\Big(kt-\frac{\beta\pi}{2}\Big).
\end{equation}

Далі скористаємось наступним твердженням роботи \cite{S_S}.

 \textbf{Лема 1.} \emph{Нехай $\gamma\in\mathbb{R}$, а
${\lambda(k)}, \ {k=1,2,...}$ --- деяка послідовність дійсних чисел.
Тоді для довільних $N,M\in\mathbb{N} \ (N<M)$ має місце рівність}
$$
\frac{1}{M-N}\sum\limits_{k=N}^{M-1}
\sum\limits_{j=1}^{k}\lambda(j)\cos\left(jt+\gamma\right)=\sum\limits_{k=1}^{N}\lambda(k)\cos\left(kt+\gamma\right)+
$$
\begin{equation}\label{statement1}
+
\frac{1}{M-N}\sum\limits_{k=N+1}^{M-1}\left(M-k\right)\lambda(k)\cos\left(kt+\gamma\right).
\end{equation}

 Поклавши в умовах леми 1 $M=n, \
N=2n-[\eta(n)]$, $\lambda(k)\equiv1, \ \gamma=-\frac{\beta\pi}{2}$,
із (\ref{statement1}) отримаємо
\begin{equation}\label{for3}
\frac{1}{[\eta(n)]-n}\sum\limits_{k=2n-[\eta(n)]+1}^{n-1}(n-k)\cos\Big(kt-\frac{\beta\pi}{2}\Big)=
\frac{1}{[\eta(n)]-n}\sum\limits_{k=2n-[\eta(n)]}^{n-1}D_{k,\beta}(t)-D_{2n-[\eta(n)],\beta}(t).
\end{equation}
Із (\ref{q}) i (\ref{for3}) випливає рівність
\begin{equation}\label{x}
\Psi^{*}_{\beta,n}(t)=\psi(n)D_{n-1,\beta}(t)-\frac{\psi(n)}{[\eta(n)]-n}\sum\limits_{k=2n-[\eta(n)]}^{n-1}D_{k,\beta}(t)+
\sum\limits_{k=n}^{\infty}\psi(k)\cos\Big(kt-\frac{\beta\pi}{2}\Big).
\end{equation}

Застосовуючи до останнього доданку з правої частини
рівності (\ref{x}) перетворення Абеля, одержимо, що при довільному $n\in~ \mathbb{N}$
\begin{equation}\label{ad}
\sum\limits_{k=n}^{\infty}\psi(k)\cos\Big(kt-\frac{\beta\pi}{2}\Big)=\sum\limits_{k=n}^{\infty}\Delta\psi(k)
D_{k,\beta}(t)- \psi(n)D_{n-1,\beta}(t),
\end{equation}
де $\Delta\psi(k)=\psi(k)-\psi(k+1)$.

Згідно з  (\ref{x}) і (\ref{ad})
\begin{equation}\label{for1}
\Psi^{*}_{\beta,n}(t)=\sum\limits_{k=n}^{\infty}\Delta\psi(k)
D_{k,\beta}(t)-\frac{\psi(n)}{[\eta(n)]-n}\sum\limits_{k=2n-[\eta(n)]}^{n-1}D_{k,\beta}(t).
\end{equation}

В силу відомої формули
\begin{equation}\label{grad}
\sum\limits_{k=0}^{N-1}\sin(\gamma+kt)=\sin\Big(\gamma+\frac{N-1}{2}t\Big)\sin
\frac{Ny}{2}\cosec \frac{t}{2}
\end{equation}
 (див., наприклад, \cite[с. 43]{Gradshteyn}), при $N=k+1$, $\gamma=(1-\beta)\frac{\pi}{2}$, маємо
$$
D_{k,\beta}(t)=\sum\limits_{j=0}^{k}\cos\Big(kt+\frac{\beta\pi}{2}\Big)-
\frac{1}{2}\cos{\beta\pi}{2}=\frac{\cos\big(\frac{kt}{2}-
\frac{\beta\pi}{2}\big)\sin\frac{k+1}{2}t}{\sin\frac{t}{2}}-\frac{1}{2}\cos\frac{\beta\pi}{2}=
$$
\begin{equation}\label{Dkb}
=\frac{\sin\big((k+\frac{1}{2})t-
\frac{\beta\pi}{2}\big)+\cos\frac{t}{2}
\sin\frac{\beta\pi}{2}}{2\sin\frac{t}{2}}, \ \ \ 0<|t|\leq \pi.
\end{equation}

З (\ref{Dkb}) і  (\ref{for1})
одержуємо
$$
\Psi^{*}_{\beta,n}(t)=\sum\limits_{k=n}^{\infty}\Delta\psi(k)\frac{\sin\big((k+\frac{1}{2})t-
\frac{\beta\pi}{2}\big)+\cos\frac{t}{2}
\sin\frac{\beta\pi}{2}}{2\sin\frac{t}{2}}-
$$
$$
-\frac{\psi(n)}{[\eta(n)]-n}\sum\limits_{k=2n-[\eta(n)]}^{n-1}
\frac{\sin\big((k+\frac{1}{2})t-
\frac{\beta\pi}{2}\big)+\cos\frac{t}{2}
\sin\frac{\beta\pi}{2}}{2\sin\frac{t}{2}}=
$$
\begin{equation}\label{k1}
=\frac{1}{2\sin\frac{t}{2}}\Big(\sum\limits_{k=n}^{\infty}\Delta\psi(k)\sin\Big(\!\Big(k+\frac{1}{2}\Big)t-
\frac{\beta\pi}{2}\Big)-\!\frac{\psi(n)}{[\eta(n)]-n}\sum\limits_{k=2n-[\eta(n)]}^{n-1}
\sin\Big(\Big(k+\frac{1}{2}\Big)t-
\frac{\beta\pi}{2}\Big)\!\Big).
\end{equation}

Застосуємо  перетворення Абеля до першої суми з правої
частини  рівності (\ref{k1}), внаслідок чого запишемо
$$
\Psi^{*}_{\beta,n}(t)\!=\!\frac{1}{2\sin\frac{t}{2}}
\Big(\sum\limits_{k=n}^{\infty}\Delta^{2}\psi(k)\sum\limits_{j=0}^{k}\sin\Big(\Big(j+\frac{1}{2}\Big)t-
\frac{\beta\pi}{2}\Big)\!-\!\Delta\psi(n)\sum\limits_{j=0}^{n-1}\sin\Big(\Big(j+\frac{1}{2}\Big)t-
\frac{\beta\pi}{2}\Big)\!-
$$
\begin{equation}\label{for2}
-\frac{\psi(n)}{[\eta(n)]-n}\sum\limits_{k=2n-[\eta(n)]}^{n-1}
\sin\Big(\Big(k+\frac{1}{2}\Big)t-
\frac{\beta\pi}{2}\Big)\Big), \ 0<|t|\leq \pi,
\end{equation}
де $\Delta^{2}\psi(k)=\Delta\psi(k)-\Delta\psi(k+1)=\psi(k)-2\psi(k+1)+\psi(k+2)$.

\noindent Оскільки
\begin{equation}\label{ineq}
\sin \frac{t}{2}\geq\frac{t}{\pi}, \ \
 0\leq t\leq\pi,
\end{equation}
то в силу  (\ref{grad})
\begin{equation}\label{ner}
\Big|\sum\limits_{j=0}^{k}\sin\Big(\Big(j+\frac{1}{2}\Big)t-
\frac{\beta\pi}{2}\Big)\Big|\leq \frac{\pi}{|t|}, \ \ 0<|t|\leq
\pi.
\end{equation}

 З  (\ref{for2})--(\ref{ner}) маємо
\begin{equation}\label{z7}
\left|\Psi^{*}_{\beta,n}(t)\right|\leq
\frac{\pi^{2}}{2t^{2}}\Big(\sum\limits_{k=n}^{\infty}\Delta^{2}\psi(k)+
\Delta\psi(n)+\frac{2\psi(n)}{[\eta(n)]-n}\Big)=
\frac{\pi^{2}}{t^{2}}\Big(\Delta\psi(n)+\frac{\psi(n)}{[\eta(n)]-n}\Big),
\ \ 0<|t|\leq \pi.
\end{equation}
Оскільки $\psi\in\mathfrak{M}$, то
\begin{equation}\label{m5}
\Delta\psi(n)\leq
|\psi'(n)|, \ \psi'(n):=\psi'(n+0).
\end{equation}

Оцінимо значення величини $|\psi'(n)|$. Для цього нам знадобиться наступне твердження.

 \textbf{Лема 2.} \emph{Нехай $\psi\in\mathfrak{M}_{\infty}^{+}$, $\mu(t)\geq b>0$. Тоді }
\begin{equation}\label{lemma}
\frac{1}{2}\frac{b^{2}}{(b+1)^{2}}\left(\eta(t)-t\right)\leq
\frac{\psi(t)}{|\psi'(t)|}\leq4\big(1+\frac{1}{b}\big)\left(\eta(t)-t\right),
\ \ t\geq 1.
\end{equation}

\textbf{Доведення леми.} Згідно з  \cite[с. 165]{Stepanets1}
справедливе співвідношення
\begin{equation}\label{ff}
|\psi'(\eta(\eta(t)))|\left(\eta(\eta(t))-\eta(t)\right)\leq
\frac{1}{4}\psi(t)\leq|\psi'(\eta(t))|\left(\eta(\eta(t))-\eta(t)\right),
\ \psi\in\mathfrak{M}.
\end{equation}

З правої частини (\ref{ff}) випливає, що
\begin{equation}\label{fff1}
\frac{\psi(t)}{|\psi'(t)|}\leq\frac{4|\psi'(\eta(t))|}{|\psi'(t)|}\left(\eta(\eta(t))-\eta(t)\right)\leq
4\left(\eta(\eta(t))-\eta(t)\right).
\end{equation}

Згідно з формулою (11) роботи \cite{S_S}
 для довільної функції
$\psi\in\mathfrak{M}_{\infty}^{+}$   при
$\mu(t)\geq b>0$ має місце співвідношення
\begin{equation}\label{i1}
\frac{1}{2}\left(\eta(t)-t\right)\leq
\eta(\eta(t))-\eta(t)<\Big(1+\frac{1}{b}\Big)\left(\eta(t)-t\right),
\end{equation}
тому з (\ref{fff1}) отримуємо
\begin{equation}\label{ff1}
\frac{\psi(t)}{|\psi'(t)|} \leq
4\Big(1+\frac{1}{b}\Big)\left(\eta(t)-t\right).
\end{equation}

З іншого боку, в силу (\ref{ff}) і (\ref{i1})
\begin{equation}\label{ff2}
\frac{\psi(t)}{|\psi'(t)|}\geq
\frac{4|\psi'(\eta(\eta(t)))|}{|\psi'(t)|}\left(\eta(\eta(t))-\eta(t)\right)\geq
\frac{2|\psi'(\eta(\eta(t)))|}{|\psi'(t)|}\left(\eta(t)-t\right).
\end{equation}
Оскільки
$$
\psi'(t)=(4\psi(\eta(\eta(t))))'=4\psi'\left(\eta(\eta(t))\right)\eta'(\eta(t))\eta'(t), \ \
\ \psi\in\mathfrak{M},
$$
(тут і надалі $\eta'(t)=\eta'(t+0)$), то
\begin{equation}\label{k12}
\frac{\psi'\left(\eta(\eta(t))\right)}{\psi'(t)}=\frac{1}{4\eta'(\eta(t))\eta'(t)}.
\end{equation}
 Згідно з формулою (13) роботи \cite{S_S} для довільної
$\psi\in\mathfrak{M}_{\infty}^{+}$ при $\mu(t)\geq b>0$
\begin{equation}\label{k14}
\eta'(t)\leq 1+\frac{1}{b}, \ t\geq 1.
\end{equation}
Тому з (\ref{k12}) i (\ref{k14}) випливає
\begin{equation}\label{m2}
\frac{2|\psi'(\eta(\eta(t)))|}{|\psi'(t)|}\left(\eta(t)-t\right)\geq
\frac{1}{2}\frac{b^{2}}{(b+1)^{2}}\left(\eta(t)-t\right).
\end{equation}

Об'єднавши (\ref{ff2}) i (\ref{m2}), отримуємо твердження леми. Лему
доведено.

В силу (\ref{m5}) i (\ref{lemma}) маємо
\begin{equation}\label{m3}
\Delta\psi(n)\leq|\psi'(n)|\leq \frac{2(b+1)^{2}}{b^{2}}\frac{\psi(n)}{\eta(n)-n}, \
\psi\in\mathfrak{M}_{\infty}^{+}, \ b>0.
\end{equation}
Згідно з лемою 2 роботи \cite{S_S}  якщо
 $\psi\in\mathfrak{M}_{\infty}^{+}$, $\eta(n)-n\geq a>1$, $\mu(n)\geq b>0$, то
\begin{equation}\label{i2}
\Big(1-\frac{1}{a}\Big)\left(\eta(n)-n\right)<[\eta(n)]-n.
\end{equation}
 З (\ref{z7}), (\ref{m3}), i (\ref{i2})
випливає нерівність
\begin{equation}\label{z8}
|\Psi^{*}_{\beta,n}(t)|\leq
\pi^{2}\Big(\frac{2(b+1)^{2}}{b^{2}}+\frac{a}{a-1}\Big)\frac{\psi(n)}{\eta(n)-n}\frac{1}{t^{2}},
\  0< t\leq\pi.
\end{equation}
Покажемо також,  що при $\eta(n)-n\geq a>0$ і $\mu(n)\geq b>2$ для довільних $t\in\mathbb{R}$ виконується
нерівність
\begin{equation}\label{z9}
|\Psi^{*}_{\beta,n}(t)|\leq
\Big(\frac{2b}{b-2}+\frac{1}{a}+\frac{1}{2}\Big)\psi(n)(\eta(n)-n).
\end{equation}
З (\ref{kernel}) маємо
\begin{equation}\label{m6}
|\Psi^{*}_{\beta,n}(t)|\leq\Big|\sum\limits_{k=n}^{\infty}\psi(k)\cos\Big(kt-\frac{\beta\pi}{2}\Big)\Big|+
\psi(n)\sum\limits_{k=2n-[\eta(n)]+1}^{n-1}\Big(1-\frac{n-k}{[\eta(n)]-n}\Big).
\end{equation}
Відповідно до формули (30) роботи \cite{S_S}  для довільних
$\psi\in\mathfrak{M}_{\infty}^{+}, \ \beta\in\mathbb{R} $ при ${\eta(n)-n\geq a>0}$ і $\mu(n)\geq b>2$
\begin{equation}\label{z10}
\Big|\sum\limits_{k=n}^{\infty}\psi(k)\cos\Big(kt-\frac{\beta\pi}{2}\Big)\Big|\leq
\Big(\frac{2b}{b-2}+\frac{1}{a}\Big)\psi(n)(\eta(n)-n), \ \
t\in\mathbb{R}.
\end{equation}
Крім того, як неважко переконатись,
$$
\psi(n)\sum\limits_{k=2n-[\eta(n)]+1}^{n-1}\Big(1-\frac{n-k}{[\eta(n)]-n}\Big)=
\psi(n)\Big([\eta(n)]-n-1-\frac{([\eta(n)]-n)([\eta(n)]-n-1)}{2([\eta(n)]-n)}\Big)=
$$
\begin{equation}\label{m7}
=\frac{\psi(n)}{2}\left([\eta(n)]-n-1\right)<\frac{\psi(n)}{2}(\eta(n)-n).
\end{equation}

З (\ref{m6})---(\ref{m7}), отримуємо нерівність (\ref{z9}).

Враховуючи (\ref{z8}), (\ref{z9}) а також нерівності
$$
\frac{a}{a-1}>\frac{1}{a}+\frac{1}{2}, \ \ a>1, \ \
\frac{b}{b-2}>\frac{(b+1)^{2}}{b^{2}}, \ \ b>2,
$$
маємо, що при $1\leq p'<\infty$, $\eta(n)-n\geq a>1$ і $\mu(n)\geq b>2$
$$
\|\Psi^{*}_{\beta,n}(t)\|_{p'} \leq
\psi(n)\Bigg(\Big(\frac{2b}{b-2}+\frac{1}{a}+\frac{1}{2}\Big)^{p'}\int\limits_{|t|\leq
\frac{1}{\eta(n)-n}}(\eta(n)-n)^{p'}dt+
$$
$$
+\pi^{2p'}\Big(\frac{2(b+1)^{2}}{b^{2}}+\frac{a}{a-1}\Big)^{p'}\frac{1}{(\eta(n)-n)^{p'}}\int\limits_{
\frac{1}{\eta(n)-n}\leq |t|\leq
\pi}\frac{dt}{t^{2p'}}\Bigg)^{\frac{1}{p'}}<
$$
$$
<\psi(n)\Big(\frac{2b}{b-2}+\frac{a}{a-1}\Big)\Bigg(\int\limits_{|t|\leq
\frac{1}{\eta(n)-n}}(\eta(n)-n)^{p'}dt+\frac{\pi^{2p'}}{(\eta(n)-n)^{p'}}\int\limits_{
\frac{1}{\eta(n)-n}\leq |t|\leq
\pi}\frac{dt}{t^{2p'}}\Bigg)^{\frac{1}{p'}}<
$$
$$
<\psi(n)(\eta(n)-n)^{1-\frac{1}{p'}}\Big(\frac{2b}{b-2}+\frac{a}{a-1}\Big)2^{\frac{1}{p'}}
\Big(1+\pi^{2p'}\frac{1}{2p'-1}\Big)^{\frac{1}{p'}}\leq
$$
\begin{equation}\label{z12}
\leq
2(1+\pi^{2})\Big(\frac{2b}{b-2}+\frac{a}{a-1}\Big)\psi(n)(\eta(n)-n)^{1-\frac{1}{p'}}.
\end{equation}

У випадку $p'=\infty$ зі співвідношення (\ref{z9}) маємо, що при  $\eta(n)-n\geq a>1$ і $\mu(n)\geq b>2$
$$
\|\Psi^{*}_{\beta,n}(t)\|_{p'}=\|\Psi^{*}_{\beta,n}(t)\|_{\infty}\leq
\Big(\frac{2b}{b-2}+\frac{1}{a}+\frac{1}{2}\Big)\psi(n)(\eta(n)-n)<
$$
\begin{equation}\label{z13}
<2(1+\pi^{2})\Big(\frac{2b}{b-2}+\frac{a}{a-1}\Big)\psi(n)(\eta(n)-n).
\end{equation}

Зі співвідношень (\ref{fff}), (\ref{zz2}), (\ref{z12}) i (\ref{z13})
випливає справедливість нерівності (\ref{dop1}).
 Теорему~1 доведено.

\textbf{Наслідок 1.} \emph{ Нехай
$\psi\in\mathfrak{M}^{+}_{\infty}$,
$\lim\limits_{n\rightarrow\infty}(\eta(\psi,n)-n)=\infty$, $\beta\in
\mathbb{R}$, ${1\leq p\leq \infty}$. Тоді
$$
{
E}_{n}(C^{\psi}_{\beta,p})_{C}\asymp
\mathcal{E}_{n}(C^{\psi}_{\beta,p})_{C}\asymp\psi(n)(\eta(n)-n)^{\frac{1}{p}}.
$$
}
Неважко переконатися, що для функції
${\psi_{r,\alpha}(t)=\exp(-\alpha t^{r})}, \ {\alpha>0}, \
{r\in(0,1)}$
\begin{equation}\label{eta-n}
\eta(n)-n=\eta(\psi_{r,\alpha};n)-n=
 n\Big(\Big(1+\frac{\ln 2}{\alpha n^{r}}\Big)^{\frac{1}{r}}-1\Big),
\end{equation}
$$\mu(n)=\mu(\psi_{r,\alpha};n)=\frac{n}{\eta(\psi_{r,\alpha};n)-n}=\frac{1}{\left(\frac{\ln 2}{\alpha n^{r}}+1\right)^{\frac{1}{r}}-1},$$
і, як показано в \cite[с. 257--258]{S_S}, при
$$
n\geq\max\Big\{1+\Big(\frac{2r\alpha}{\ln2}\Big)^{\frac{1}{1-r}},
1+2\Big(\frac{\ln2}{\alpha(3^{r}-2^{r})}\Big)^{\frac{1}{r}}\Big\}
$$
виконуються нерівності
$$
\eta(n)-n\geq a(\alpha,r)>2,
$$
$$
\mu(n)\geq a(\alpha,r)>2,
$$
де
\begin{equation}\label{aa}
a(\alpha,r)=\frac{\ln 2}{\alpha
r}\Big(1+\Big(\frac{2r\alpha}{\ln2}\Big)^{\frac{1}{1-r}}\Big)^{1-r},
\end{equation}
\begin{equation}\label{bb}
b(\alpha,r)=\Big(\Big(\frac{\ln 2}{\alpha
}\Big(1+2\Big(\frac{\ln2}{\alpha(3^{r}-2^{r})}\Big)
^{\frac{1}{r}}\Big)^{-r}+1\Big)^{\frac{1}{r}}-1\Big)^{-1}.
\end{equation}
Тоді з теореми 1 випливає наступне твердження.

 \textbf{Наслідок 2.} \emph{ Нехай ${r\in(0,1)}, \ {\alpha>0}$, ${1\leq  p\leq
\infty}$ і $\beta\in \mathbb{R}$. Тоді  для всіх $n$ таких, що
\begin{equation}\label{n}
n\geq\max\Big\{1+\Big(\frac{2r\alpha}{\ln2}\Big)^{\frac{1}{1-r}},
1+2\Big(\frac{\ln2}{\alpha(3^{r}-2^{r})}\Big)^{\frac{1}{r}}\Big\},
\end{equation}
 справедливі оцінки}
\begin{equation}\label{cons}
C_{a} \exp\left(-\alpha
n^{r}\right)n^{\frac{1}{p}}\Big(\!\!\Big(1+\frac{\ln 2}{\alpha
n^{r}}\Big)^{\frac{1}{r}}\!-\!1\!\Big)^{\frac{1}{p}} \leq {
E}_{n}\Big(C^{\alpha,r}_{\beta,p}\Big)_{C}\leq C^{*}_{a,b}
\exp\left(-\alpha
n^{r}\right)n^{\frac{1}{p}}\Big(\!\!\Big(1+\frac{\ln 2}{\alpha
n^{r}}\Big)^{\frac{1}{r}}\!-\!1\!\Big)^{\frac{1}{p}},
 \end{equation}
\emph {де величини $C_{a}$ i $C^{*}_{a,b}$
означаються формулами (\ref{Ca}) i (\ref{Cab}) при
${a=a(\alpha,r)}$, ${b=b(\alpha,r)}$, що задані за допомогою
рівностей (\ref{aa}) і (\ref{bb}).}

\textbf{Теорема 2.} \emph{ Нехай $\psi\in\mathfrak{M}^{+}_{\infty}$,
$\beta\in \mathbb{R}$, ${1\leq s\leq\infty}, \
\frac{1}{s}+\frac{1}{s'}=1$. Тоді  для довільних $n\in
 \mathbb{N}$, таких, що $ {\eta(n)-n\geq a>2}, \ {\mu(n)\geq b>2}$ справедливі оцінки
\begin{equation}\label{theorem2}
C_{a}\psi(n)(\eta(n)-n)^{\frac{1}{s'}}\leq{
E}_{n}(L^{\psi}_{\beta,1})_{s}\leq C^{*}_{a,b}
\psi(n)(\eta(n)-n)^{\frac{1}{s'}},
\end{equation}
де сталі $C_{a}$ i $C^{*}_{a,b}$ означаються
формулами (\ref{Ca}) i (\ref{Cab}) відповідно.}

\textbf{\emph{Доведення теореми 2.}} Згідно з формулою (70) роботи \cite{S_S} для довільних
${\psi\in\mathfrak{M}^{+}_{\infty}}$ при $ {\eta(n)-n\geq a>2}, \ {\mu(n)\geq b>2}$ має місце оцінка
$$
{
E}_{n}(L^{\psi}_{\beta,1})_{s}\geq C_{a}\psi(n)(\eta(n)-n)^{\frac{1}{s'}},
$$
в якій величина $C_{a}$ означена рівністю (\ref{Ca}). Для оцінки зверху найкращих наближень
${
E}_{n}\big(L^{\psi}_{\beta,1}\big)_{s}$ розглянемо величину
$$
{\cal
E}(L^{\psi}_{\beta,1};V_{n,\psi})_{s}=\sup\limits_{f\in
L^{\psi}_{\beta,1}}\|f(\cdot)-V_{n,\psi}(f;\cdot)\|_{s}, \ \ 1\leq
s\leq \infty,
$$
де суми $V_{n,\psi}$ означаються формулою (\ref{sum}).

\noindent Оскільки
\begin{equation}\label{dop}
{ E}_{n}(L^{\psi}_{\beta,1})_{s}\leq {\cal
E}(L^{\psi}_{\beta,1};V_{n,\psi})_{s},
\end{equation}
то для доведення теореми 2 достатньо показати справедливість співвідношення
\begin{equation}\label{sss}
{\cal
E}(L^{\psi}_{\beta,1};V_{n,\psi})_{s}
\leq C_{a,b}^{*}
\psi(n)(\eta(n)-n)^{\frac{1}{s'}}.
\end{equation}
\noindent
Для цього використаємо інтегральне зображення (\ref{fff}), яке у випадку $f\in L^{\psi}_{\beta,1}$ буде справедливим майже для всіх $x\in\mathbb{R}$,  та
нерівність (5.28) з \cite[с. 43]{Korn}. Тоді для довільних $1\leq
s\leq\infty$, одержимо
\begin{equation}\label{t1}
{ E}_{n}\big(L^{\psi}_{\beta,1}\big)_{s}\leq
\frac{1}{\pi}\|\Psi^{*}_{\beta,n}(\cdot)\|_{s}\|
\varphi(\cdot)\|_{1}\leq\frac{1}{\pi}\|\Psi^{*}_{\beta,n}(\cdot)\|_{s}.
\end{equation}
Скориставшись співвідношеннями (\ref{z12}) i (\ref{z13}), поклавши $p'=s$, з (\ref{t1}) отримуємо
(\ref{sss}).
Теорему 2 доведено.

\textbf{Наслідок 3.} \emph{ Нехай
$\psi\in\mathfrak{M}^{+}_{\infty}$,
$\lim\limits_{n\rightarrow\infty}(\eta(\psi,n)-n)=\infty$, $\beta\in
\mathbb{R}$, ${1\leq  s\leq \infty}$, $\frac{1}{s}+\frac{1}{s'}=1$.
Тоді  }
$$
{
E}_{n}(L^{\psi}_{\beta,1})_{s}\asymp \mathcal{E}_{n}(L^{\psi}_{\beta,1})_{s}\asymp\psi(n)(\eta(n)-n)^{\frac{1}{s'}}.
$$

\textbf{Наслідок 4.} \emph{ Нехай ${r\in(0,1)}, \ {\alpha>0}$, ${1\leq  s\leq
\infty}$, $\frac{1}{s}+\frac{1}{s'}=1$, $\beta\in \mathbb{R}$. Тоді
для всіх  $n\in\mathbb{N}$ таких, що задовольняють умову (\ref{n}),
 справедливі оцінки}
$$
C_{a} \exp\left(-\alpha
n^{r}\right)n^{\frac{1}{s'}}\Big(\!\!\Big(1+\frac{\ln 2}{\alpha
n^{r}}\Big)^{\frac{1}{r}}\!-\!1\!\Big)^{\frac{1}{s'}} \leq {
E}_{n}(L^{\alpha,r}_{\beta,1})_{s}\leq C^{*}_{a,b}
\exp\left(-\alpha
n^{r}\right)n^{\frac{1}{s'}}\Big(\!\!\Big(1+\frac{\ln 2}{\alpha
n^{r}}\Big)^{\frac{1}{r}}\!-\!1\!\Big)^{\frac{1}{s'}},
$$
\emph {де величини $C_{a}$ i $C^{*}_{a,b}$
означаються формулами (\ref{Ca}) i (\ref{Cab}) при
${a=a(\alpha,r)}$, $b=b(\alpha,r)$, які задані  рівностями (\ref{aa}) i (\ref{bb}).}

З рівності (\ref{eta-n}) для величини $\eta(\psi_{r,\alpha};n)-n$ неважко одержати двосторонні оцінки
\begin{equation}\label{q123}
\frac{\ln2}{\alpha r}n^{1-r}\leq\eta(\psi_{r,\alpha};n)-n\leq (1+\ln2^{1/\alpha})^{\frac{1-r}{r}}
\frac{\ln2}{\alpha r}n^{1-r}, \ \alpha>0, \
r\in(0,1), \ \ n\in\mathbb{N}.
\end{equation}

З  наслідків 2 і 4 та формули (\ref{q123}) випливають
 порядкові рівності
$$
{ E}_{n}(C^{\alpha,r}_{\beta,p})_{C}\asymp\exp\left(-\alpha
n^{r}\right)n^{\frac{1-r}{p}},\  \ \ \ 1\leq p\leq  \infty,
$$
$$
{ E}_{n}(L^{\alpha,r}_{\beta,1})_{s}\asymp\exp\left(-\alpha
n^{r}\right)n^{\frac{1-r}{s'}}, \  \ \ {1\leq  s\leq \infty},\
{\frac{1}{s}+\frac{1}{s'}=1}.
$$

 Співставляючи співвідношення (\ref{th1}) i
(\ref{th2}) з (\ref{theorem1}) i (\ref{theorem2}) відповідно, при
виконанні умов теорем 1 та 2 можна записати оцінки
$$
C_{a}\psi(n)(\eta(n)-n)^{\frac{1}{p}}\leq{
E}_{n}(C^{\psi}_{\beta,p})_{C}\leq
C_{a,b}(p)\psi(n)(\eta(n)-n)^{\frac{1}{p}},
$$
$$
C_{a}\psi(n)(\eta(n)-n)^{\frac{1}{s'}}\leq{
E}_{n}(L^{\psi}_{\beta,1})_{s}\leq C_{a,b}(s')
\psi(n)(\eta(n)-n)^{\frac{1}{s'}}, \ \
\frac{1}{s}+\frac{1}{s'}=1,
$$
де
$$
C_{a,b}(p)=\min\big\{(2p)^{1-\frac{1}{p}} \ C_{a,b}, \
C^{*}_{a,b}\big\},
$$
а величини $C_{a,b}$ i $C_{a,b}^{*}$ означені формулами (\ref{Cab_F}) i (\ref{Cab}) відповідно.

Обчислення показують, що при невеликих значеннях $p \ (1\leq
p\leq 7)$, $ \ \ {C_{a,b,p}=(2p)^{1-\frac{1}{p}}} \ C_{a,b}$, а при великих $p \ (p\geq 26)$, $ \  \ {C_{a,b,p}=C^{*}_{a,b}}$.

E-mail: \href{mailto:serdyuk@imath.kiev.ua}{serdyuk@imath.kiev.ua},
\href{mailto:tania_stepaniuk@ukr.net}{tania$_{-}$stepaniuk@ukr.net}

\end{document}